\documentclass[a4paper,twoside]{amsart}
\usepackage{amssymb}
\usepackage{amscd}
\usepackage{eucal}
\usepackage{amsfonts}
\usepackage{latexsym}
\usepackage{amsthm}
\newtheorem{tw}{Theorem}[section]

\newtheorem{wn}[tw]{Corollary}
\newtheorem{fakt}[tw]{Fact}
\newtheorem{df}[tw]{Definition}

\newtheorem{problem}[tw]{Problem}

\newcommand{\lif}{\liminf_{\sigma \in \Sigma}}

\newcommand{\e}{\varepsilon}
\newcommand{\f}{\varphi}

\newcommand{\D}{\Delta}

\newcommand{\N}{\mathbb{N}}

\begin{document}
\title{Coarse version of the Banach-Stone theorem}
\author{Rafa\l ~G\'orak}
\address{Rafa\l ~G\'orak\\ Technical University of Warsaw \\
Pl. Politechniki 1\\  00-661 Warszawa \\ Poland \\ email:
rafalgorak@gmail.com}
\date{\today}
\subjclass[2000]{46E15, 46B26, 46T99} \keywords{Banach-Stone
theorem, function space, isometry}

\begin{abstract}
We show that if there exists a Lipschitz homeomorphism $T$ between
the nets in the Banach spaces $C(X)$ and $C(Y)$ of continuous real
valued functions on compact spaces $X$ and $Y$, then the spaces
$X$ and $Y$ are homeomorphic provided $l(T) \times
l(T^{-1})<\frac{6}{5}$. By $l(T)$ and $l(T^{-1})$ we denote the
Lipschitz constants of the maps $T$ and $T^{-1}$. This improves
the classical result of Jarosz and the recent result of Dutrieux
and Kalton where the constant obtained is $\frac{17}{16}$. We also
estimate the distance of the map $T$ from the isometry of the
spaces $C(X)$ and $C(Y)$.
\end{abstract}

\maketitle

\section{Introduction}
This paper deals with Banach spaces and some equivalences arising
from the concepts of large scale geometry (coarse geometry) and
the geometry of Banach spaces. Let us first recall the very well
known notion of a Lipschitz map and the related notion of a
bi-Lipschitz map:
\begin{df}\label{lip}
Let $E$ and $F$ be metric spaces and $T:E \mapsto F$ be a map
between these spaces. Then $T$ is
\begin{itemize}
\item a Lipschitz map if for some constant $M>0$ it satisfies
 the inequality $d_F(Tx,Ty) \leq Md_E(x,y)$ for all $x,y \in E$; \item
 an
$M$ - bi-Lipschitz map if it is an onto map and the inequality
$$\frac{1}{M}d_E(x,y) \leq d_F(Tx,Ty) \leq Md_E(x,y)$$ holds for
all $x,y \in E$.
\end{itemize}
For a Lipschitz map $T$ we denote by $l(T)= \sup \{
\frac{d_F(Tx,Ty)}{d_E(x,y)}; \; x \neq y \}$.
\end{df}
Let us now introduce two notions that are generalizations of the
above in the large scale direction:

\begin{df}\label{coarse lip}
Let $E$ and $F$ be metric spaces and $T:E \mapsto F$ be a map
between these spaces. Then $T$ is a coarse Lipschitz map
(Lipschitz for large distances) if for some constants $M>0$ and $L
\geq 0$
 the inequality \mbox{$d_F(Tx,Ty) \leq Md_E(x,y)+L$} is satisfied for all $x,y \in E$.
For a coarse Lipschitz map $T$ we denote by $l_{\infty}(T)$ the
number $\inf_{\theta>0} \sup_{\|x-y\| \geq \theta}
\frac{\|Tx-Ty\|}{\|x-y\|}$ (a Lipschitz constant at $\infty$).
\end{df}
It is worth mentioning that if $T$ is a uniformly continuous map
between Banach spaces then it is also a coarse Lipschitz map.
Although this is a basic fact it is very important for a nonlinear
classification of Banach spaces. For more information please see
the book \cite{BeLin}. Another definition is a natural
generalization of the notion of bi-Lipschitz map.
\begin{df}\label{def quasi}
Let $E$ and $F$ be metric spaces and $T:E \mapsto F$ be a map
between these spaces. $T$ is a coarse $(M,L)$-quasi isometry (or
just coarse quasi isometry)  if it satisfies the following
conditions:
\begin{itemize}
\item $\frac{1}{M}d_E(x,y) -L \leq d_F(Tx,Ty) \leq Md_E(x,y)+L$
for all $x,y \in E$; \item there exists $\xi>0$ such that for
every $y \in F$ there exists $x \in E$ such that $d_F(y,Tx) \leq
\xi$. In other words $T(E)$ is $\xi$ dense in $F$.
\end{itemize}

\end{df}
The class of maps defined above is sometimes called just quasi
isometries (see \cite{Bridson}) however the "coarse" is added
since in some contexts quasi isometries are those coarse quasi
isometries for which $L=0$ (see \cite{Vest1}). We shall focus on
coarse quasi isometries between Banach spaces.

In our considerations we find a condition under which the
existence of a coarse quasi isometry between spaces $C(X)$ and
$C(Y)$ implies the existence of a homeomorphism of topological
spaces $X$ and $Y$. The last section of this paper is devoted to
stability problems i.e. we estimate the distance of a coarse
$(M,L)$-quasi isometry of the Banach spaces $C(X)$ and $C(Y)$ to
an isometry of these spaces as $M \rightarrow 1$. In paper
\cite{DutKal} Dutrieux and Kalton consider different kinds of
nonlinear distances between Banach spaces. Among the others they
consider the uniform distance between two Banach spaces $E$ and
$F$ as well as the net distance. Let us recall the definition of
both:
\begin{df}
Let $E$ and $F$ be Banach spaces. By $d_u(E,F)$, $d_N(E,F)$ we
denote the uniform and the net distance, respectively.
\begin{itemize}
\item $d_u(E,F) = \inf l_{\infty}(u) \times l_{\infty}(u^{-1})$,
where the infimum is taken over all uniform homeomorphisms $u$
between $E$ and $F$. \item $d_N(E,F) = \inf l(T) \times
l(T^{-1})$, where the infimum is taken over all bi-Lipschitz maps
$T$ between the nets $N_E$ and $N_F$ in the Banach spaces $E$,
$F$, respectively.
\end{itemize}
\end{df}

Let us recall that a subset $N_E \subset E$ of a metric space $E$
is called an $(\e,\delta)$ - net if every element of $E$ is of a
distance less then $\e$ to some element of $N_E$. Moreover every
two elements of $N_E$ are of the distance at least $\delta$. It is
easy to observe that $d_u(E,F) \geq d_N(E,F)$. In the mentioned
paper of Dutrieux and Kalton they work with the Gromov-Hausdorff
distance $d_{GH}(E,F)$ between Banach spaces $E$ and $F$. They
show that $d_N(E,F)\geq d_{GH}(E,F)+1$. This fact and their result
that the inequality $d_{GH}(C(X),C(Y)) <\frac{1}{16}$ implies the
existence of a homeomorphism of compact spaces $X$ and $Y$ give us
that $X$ and $Y$ are also homeomorphic when $d_{N}(C(X),C(Y))
<\frac{17}{16}$. In our paper we improve the constant to
$\frac{6}{5}$. Please keep in mind that we consider only the the
net distance.

Let us now discuss the connection between the net distance and the
notion of a coarse quasi isometry. We shall start with the
following fact:
\begin{fakt}\label{fakt o przedluzaniu sieci}
 Let us consider a coarse $(M,L)$-quasi isometry $T: A \mapsto B$ from a $\xi_E$ dense set in Banach spaces $E$
 onto
 a
  $\xi_F$ dense set in $F$. Then there exists a map $\widetilde{T}:E \mapsto F$, which is a bijective coarse $(M,(4M^2+3)L+4\xi_F+2M\xi_E)$-quasi isometry and
  $\|\widetilde{T}x-Tx\| \leq (2M^2+2)L+2\xi_F+M\xi_E$ for all $x \in A$.
\end{fakt}
\begin{proof}
Consider $N_E$ a maximal $ML+ M\eta$ separated set in $A$ where
$\eta>0$ is arbitrary. Let us first notice that $N_E$ is an
$(\e_E,\delta_E)$ net in E where \mbox{$\e_E=ML+ M\eta+\xi_E$} and
$\delta_E=ML+ M\eta$. Obviously $T|N_E$ is a bijection between the
net $N_E$ and $N_F=T(N_E)$ which is also an $(\e_F,\delta_F)$ net
in $F$ where $\e_F=M^2L+ M^2\eta+L+\xi_F$ (since $N_E$ is $ML+
M\eta$ dense in $A$ and $T$ is onto) and $\delta_F=\eta$. Let us
enumerate the elements of $N_E$ that is
$N_E=\{x_{\alpha}\}_{\alpha \in \tau}$. Modifying balls around the
points of $N_E$ and $N_F$  we obtain families
$(E_{\alpha})_{\alpha \in \tau}$, $(F_{\alpha})_{\alpha \in \tau}$
of subsets of $E$ and $F$, respectively such that:
\begin{itemize}
\item[(i)] $B(x_{\alpha},\frac{\delta_E}{2}) \subset E_{\alpha}$,
$B(Tx_{\alpha},\frac{\delta_E}{2}) \subset F_{\alpha}$;
\item[(ii)] $E_{\alpha} \subset B(x_{\alpha}, \e_E)$ and
$F_{\alpha} \subset B(Tx_{\alpha}, \e_F)$; \item[(iii)]
$E_{\alpha} \cap E_{\beta} =\emptyset$, $F_{\alpha} \cap F_{\beta}
=\emptyset$ for $\beta \neq \alpha$; \item[(iv)] $\bigcup_{\alpha
\in \tau} E_{\alpha} = E$ and $\bigcup_{\alpha \in \tau}
F_{\alpha} = F$.
\end{itemize}
Since $N_E$ and $N_F$ are of the same cardinality then also all
open sets in $E$ and $F$ have the same cardinalities. Hence for
every $\alpha \in \tau$ we can find a bijection
$\widetilde{T}_{\alpha}: E_{\alpha} \mapsto F_{\alpha}$ such that
$\widetilde{T}_{\alpha}(x_{\alpha})=Tx_{\alpha}$. Setting
$\widetilde{T}=\bigcup_{\alpha \in \tau} \widetilde{T}_{\alpha}$
we obtain the desired map. Indeed let us take $x, y \in E$. Choose
$E_{\alpha}$ and $E_{\beta}$ such that $x \in E_{\alpha}$ and $y
\in E_{\beta}$. We have then
$$\frac{1}{M}\|x_{\alpha}-
x_{\beta}\|-L \leq
\|\widetilde{T}x_{\alpha}-\widetilde{T}x_{\beta}\| \leq
M\|x_{\alpha}- x_{\beta}\|+L.$$ Since the partitions satisfy
conditions (i) and (ii) we obtain:
$$\|\widetilde{T}x_{\alpha}-\widetilde{T}x\| \leq \e_F,$$
$$\|\widetilde{T}x_{\beta}-\widetilde{T}y\| \leq \e_F,$$
$$\|x_{\alpha}-x\| \leq \e_E,$$
$$\|x_{\beta}-y\| \leq \e_E.$$
Combining all the above inequalities we show that $\widetilde{T}$
is a bijective coarse\\ \mbox{$(M,L+2M\e_E+2\e_F)$-quasi
isometry.} Consider $x \in A$. Obviously there exists
\mbox{$\alpha \in \tau$} such that $x \in E_{\alpha}$. From the
definition of $\widetilde{T}$ we know that \mbox{$\widetilde{T}x
\in F_{\alpha} \subset B(Tx_{\alpha},\e_F)$} and
$\widetilde{T}x_{\alpha} = Tx_{\alpha}$. This way we obtain the
inequalities:

$$\| \widetilde{T}x -  \widetilde{T}x_{\alpha} \| \leq \e_F,$$

$$ \|Tx-Tx_{\alpha}\| \leq M\e_E+L.$$

Hence $\| \widetilde{T}x -Tx\| \leq M\e_E+\e_F+L$. Setting  $\eta
= \frac{\xi_F}{2M^2}$ the proof is finished.
\end{proof}
At this moment it is worth mentioning  that from the definition of
$d_N(E,F)$ it follows that for every $M>0$ such that
$M^2>d_N(E,F)$ there are nets $N_E$ and $N_F$ in $E$ and $F$,
respectively and a bi-Lipschitz map $T$ between them such that
\mbox{$M^2 = l(T) \times l(T^{-1})
> d_N(E,F)$}. However we cannot be sure that \mbox{$l(T)=l(T^{-1})$}.
In order to obtain that consider the map $\widetilde{T}x =
\sqrt{\frac{l(T^{-1})}{l(T)}}Tx$. It is a bi-Lipschitz map of nets
(different ones) in $E$ and $F$ such that
$l(\widetilde{T})=l(\widetilde{T}^{-1})=M$. From this remark and
Fact \ref{fakt o przedluzaniu sieci} we can easily deduce
\begin{fakt}
Let $E$ and $F$ be Banach spaces. If $d_N(E,F)  < M^2$ then there
exists a bijective coarse $(M,L)$-quasi isometry between $E$ and
$F$ for some constant $L \geq 0$.
\end{fakt}
\begin{proof}
Indeed if $d_N(E,F)  < M^2$ then there exists a bijective map
$T:N_E \mapsto N_F$ between nets in Banach spaces $E$ and $F$,
respectively such that $$\frac{1}{M} \|x-y\| \leq \|Tx-Ty\| \leq
M\|x-y\|$$ for all $x,y \in N_E$. From Fact \ref{fakt o
przedluzaniu sieci} there exists $\widetilde{T}$ which is the
desired coarse $(M,L)$-quasi isometry.
\end{proof}
On the other hand it is not difficult to show that if there exists
a coarse $(M,L)$-quasi isometry between spaces $E$ and $F$ then
$d_N(E,F) \leq M^2$.

The following fact is another conclusion from Fact \ref{fakt o
przedluzaniu sieci}:
\begin{fakt}\label{fakt o bliskiej bijekcji}
Let us consider a coarse $(M,L)$-quasi isometry $T: E \mapsto F$
of Banach spaces $E$ and $F$ where $\xi$ is such that $TE$ is
$\xi$ dense in $F$. Then there exists a bijective coarse
$(M,(4M^2+3)L+4\xi)$-quasi isometry $\widetilde{T}: E \mapsto F$
such that
  $\|\widetilde{T}x-Tx\| \leq (2M^2+2)L+2\xi$ for all $x \in E$.
\end{fakt}

From now on we only consider bijective coarse $(M,L)$-quasi
isometries. Indeed from the above fact we know that whenever we
have a coarse $(M,L)$-quasi isometry we can define a new bijective
coarse quasi isometry changing only the constant $L$. Moreover it
is not more than $10L+2\xi$ from the original coarse $(M,L)$-quasi
isometry if $M<2$ which is going to be our case in further
considerations. Obviously $\xi$ is as in Definition \ref{def
quasi}. Since we are concerned mainly with large distances then
both maps, the original one and the corrected bijective version
are no different from the large scale perspective. The only thing
that matters is the constant $M$ which is unchanged.

\section{The coarse version of the Banach-Stone theorem}
As it was mentioned in Section 1 we investigate Banach spaces of
continuous real valued functions defined on compact spaces. For a
compact spaces $X$ we denote such a Banach space by $C(X)$. As
usual $C(X)$ is endowed with the $sup$ norm. Our main goal is to
consider maps between such spaces that are small perturbations of
isometries (from the large scale perspective) i.e. we consider
coarse $(M,L)$-quasi isometries where $M \rightarrow 1$.
\begin{tw}\label{tw coarse banach stone}
Let $X$ and $Y$ be compact spaces and $C(X)$, $C(Y)$ Banach spaces
of continuous real valued functions on  $X$ and $Y$, respectively.
Let $T:C(X) \mapsto C(Y)$ be a bijective coarse $(M,L)$- quasi
isometry such that $T(0)=0$. Then for every $M < \sqrt{1.2}$ there
is a homeomorphism $\varphi:X \mapsto Y$ such that for every $f
\in C(X)$, $x \in X$
$$\big{|}|Tf(\varphi(x))| - |f(x)|\big{|} \leq 5(M^2-M)\|f\| + \D.$$
The constant $\D$ depends only on $M$ and $L$. Moreover, for
$L=0$, we have $\D=0$.
\end{tw}

The proof of this theorem is a modified technique of K. Jarosz
from \cite{Jarosz}. Without loss of generality we can assume that
$L$ is so that both $T$ and $T^{-1}$ are coarse $(M,L)$-quasi
isometries. We also use the notion of Moore-Smith convergence when
dealing with topology of general topological spaces. $\Sigma$ will
always denote a directed set and whenever we write $a_{\sigma}
\rightarrow a$ we always mean $\lim_{\sigma \in \Sigma}a_{\sigma}=
a$. We give two special definitions, and then state and prove
several facts which will be used to prove the theorem.
\begin{df}
$(f^m_{\sigma})_{{\sigma} \in \Sigma} \subset C(X)$ is the m-peak
sequence at $x \in X$, for some directed set $\Sigma$ if
\begin{itemize}
\item $\|f^m_{\sigma}\|=|f^m_{\sigma}(x)|=m$ for all ${\sigma} \in
\Sigma$, \item $\lim_{{\sigma} \in \Sigma} f^m_{\sigma}|(X
\setminus U) \equiv 0$ uniformly for all open neighborhoods $U$ of
$x$.
\end{itemize}
The set of m-peak sequences at $x$ we denote by $P^X_m(x)$.
\end{df}
\begin{df} Let $D>0$ and $m>0$. We define the following:\\
$S^D_m(x)=\{y \in Y \textrm{; \;} \exists (f^m_{\sigma})_{{\sigma}
\in \Sigma} \in P^X_m(x) \textrm{\;} \exists y_{\sigma}
\rightarrow y \forall {\sigma} \in \Sigma \textrm{\;}
Tf^m_{\sigma}(y^m_{\sigma})
\geq Dm \textrm{ and }\\ T(-f^m_{\sigma})(y^m_{\sigma}) \leq -Dm\}$\\ and analogically\\
$S^{-D}_m(y)=\{x \in X \textrm{; \;} \exists
(g^m_{\sigma})_{{\sigma} \in \Sigma} \in P^Y_m(y) \textrm{\;}
\exists x^m_{\sigma} \rightarrow x \forall {\sigma} \in \Sigma
\textrm{\;} T^{-1}g^m_{\sigma}(x^m_{\sigma}) \geq Dm \textrm{ and }\\
T^{-1}(-g^m_{\sigma})(x^m_{\sigma}) \leq -Dm\}$.
\end{df}
\begin{fakt}\label{fakt ze nonempty}
Let us consider $D$ such that $D < \frac{2}{M}-M$. There exists
$m_0$ (depending on $M$, $L$ and $D$) such that for all $m > m_0$
we have $S^D_m(x) \neq \emptyset$ for all $x \in X$. Moreover if
$L=0$ then $m_0=0$.
\end{fakt}
\begin{proof}
Let us take any $(\widetilde{f}^m_{\sigma})_{\sigma \in \Sigma}
\in P^X_m(x)$ such that for all $\sigma \in \Sigma$
$\widetilde{f}^m_{\sigma}(x) = m$.  We have
$$\forall \sigma \in \Sigma \; \|T\widetilde{f}^m_{\sigma}-T(-\widetilde{f}^m_{\sigma})\| \geq
\frac{2}{M}m-L.$$ Hence $\forall \sigma \in \Sigma$ there exists
$y^m_{\sigma} \in Y$ such that
$|T\widetilde{f}^m_{\sigma}(y^m_{\sigma})-T(-\widetilde{f}^m_{\sigma})(y^m_{\sigma})|
\geq \frac{2}{M}m-L$. Let us observe that numbers
$T\widetilde{f}^m_{\sigma}(y^m_{\sigma})$ and
$T(-\widetilde{f}^m_{\sigma})(y^m_{\sigma})$ must be of different
signs. Assume the contrary. Since $\|T(\pm
\widetilde{f}^m_{\sigma})\| \leq Mm+L$ we have $Mm+L \geq
\frac{2}{M}m-L$ which is impossible for $m$ large enough provided
$\frac{2}{M}>M$ (that is if $M < \sqrt{2}$). We can and we do
assume that $\forall \sigma \in \Sigma$
$T\widetilde{f}^m_{\sigma}(y^m_{\sigma}) \geq 0$ or $\forall
\sigma \in \Sigma$ $T\widetilde{f}^m_{\sigma}(y^m_{\sigma}) \leq
0$. We define $\forall \sigma \in \Sigma \;
f^m_{\sigma}=\widetilde{f}^m_{\sigma}$ if
$T\widetilde{f}^m_{\sigma}(y^m_{\sigma}) \geq 0$ or $\forall
\sigma \in \Sigma \; f^m_{\sigma}=-\widetilde{f}_{\sigma}$
otherwise. We have
$Tf^m_{\sigma}(y^m_{\sigma})-T(-f^m_{\sigma})(y^m_{\sigma}) \geq
\frac{2}{M}m-L$. Because $\|T(\pm f^m_{\sigma})\| \leq Mm+L$ then
$$Tf^m_{\sigma}(y^m_{\sigma}) \geq (\frac{2}{M}-M)m-2L$$
$$T(-f^m_{\sigma})(y^m_{\sigma}) \leq -(\frac{2}{M}-M)m+2L.$$
By compactness of $Y$ we can assume that $y^m_{\sigma} \rightarrow
y \in Y$. Therefore for every $D<\frac{2}{M}-M$ there exists such
$m_0$(depending on $D$, $M$ and $L$) that $S^D_m(x) \neq
\emptyset$ for $m>m_0$. Let us notice that for $L=0$ we have
$m_0=0$.
\end{proof}
\begin{fakt}\label{nierownosc Tf<f} For every $m > m_0$ and $f$ such that $m \geq \|f\|$ we have
$$ |Tf(y)| \leq |f(x)| + \e(M)m+L$$
where $\e(M)=2M-1-D$ and $y$ is any element of $S^D_m(x)$.
\end{fakt}
\begin{proof}
Consider $(f^m_{\sigma})_{{\sigma} \in \Sigma} \in P^X_m(x)$ and
the corresponding sequence $y^m_{\sigma} \rightarrow y$ as in the
definition of $S^D_m(x)$. We have that $\lif
\||f^m_{\sigma}|+|f|\| = m+|f(x)|$. Since
$\max\{\|f+f^m_{\sigma}\|,\|f-f^m_{\sigma}\|\} =
\||f^m_{\sigma}|+|f|\|$ we get: $$\lif \max\{\|Tf -
T(-f^m_{\sigma})\|,\|Tf-Tf^m_{\sigma}\|\} \leq
Mm+(M-1)|f(x)|+|f(x)|+L.$$ Using the facts that
$|Tf(y^m_{\sigma})|+Dm \leq |Tf(y^m_{\sigma}) - T(\lambda
f^m_{\sigma})(y^m_{\sigma})|$ for some $\lambda \in \{-1,1\}$,
$\|f\| \leq m$ and $y^m_{\sigma} \rightarrow y$ we obtain:
\begin{eqnarray*}
Dm+|Tf(y)| & \leq &
\lif \max\{|Tf(y^m_{\sigma}) - T(-f^m_{\sigma})(y^m_{\sigma})|,|Tf(y^m_{\sigma})-Tf^m_{\sigma}(y^m_{\sigma})|\}\\
& \leq &\lif \max\{\|Tf - T(-f^m_{\sigma})\|,\|Tf-Tf^m_{\sigma}\|\}\\
& \leq & Mm+(M-1)|f(x)|+|f(x)|+L \leq (2M-1)m+|f(x)|+L.
\end{eqnarray*}
\end{proof}

\begin{fakt}\label{fi o psi = id}
Assume that $1-\e(M)M-\e(M) > 0$ where $\e(M)=2M-1-D$. Then for
every $\delta
> 0$ there exists $m_1 \geq m_0$ such that if $x_0 \in X$, $m>m_1$
and $Mm+L+ \delta>k \geq Mm+L$ then $S^{-D}_k(y)=\{x_0\}$ provided
$y \in S^D_m(x_0)$. Moreover the choice of $m_1$ depends on $M$,
$L$, $D$ and $\delta$ only and if $L=0$ then $m_1$ can be as close
to $0$ as we wish.
\end{fakt}
\begin{proof}
Let us assume the contrary that $x_1 \in S^{-D}_k(y)$ and $x_1
\neq x_0$. Consider $f$ such that $f(x_0)=0$ and $f(x_1)=m=\|f\|$.
Applying Fact \ref{nierownosc Tf<f} to $T$ and $T^{-1}$ we obtain
$$m-\e(M)k-L=|f(x_1)| - \e(M)k -L \leq |Tf(y)| \leq \e(M)m + |f(x_0)|
+L=\e(M)m+L$$ if only $k \geq Mm+L \geq \|Tf\|$. After rearranging
we have
$$m-\e(M)(m+Mm+L+\delta) \leq m-\e(M)(m+k) \leq 2L.$$
This is equivalent to the condition
$$(1-\e(M)M-\e(M))m \leq 2L+\e(M)L+\e(M)\delta.$$
Therefore  if $1-\e(M)M-\e(M) > 0$ we get a contradiction for
sufficiently large $m$. If $L=0$ then if we want to obtain $m_1$
as small as we wish it is enough to consider $\delta$ small
enough.
\end{proof}
Let us assume from now on that $M$ and $D$ are such that $D$
satisfies the condition from Fact \ref{fakt ze nonempty} and the
inequality $1-\e(M)M-\e(M) > 0$ holds. Consider \mbox{$m >
Mm_1+L$}. Since the construction of $S^{-D}_m$ and $S^D_m$ is
symmetric and involves exactly the same constants we can conclude
from the above fact that $\bigcup_{x \in X} S^D_m(x) = Y$ and
$\bigcup_{y \in Y} S^{-D}_m(y) = X$. Moreover for every $x \in X$
and $y \in Y$ we have that $S^D_m(x)$ and $S^{-D}_m(y)$ are
singletons. Hence we can clearly see that $S^D_{k}(x) =
S^D_{l}(x)$ provided $k,l \geq Mm_1+L$ and $|k-l|<\delta$. However
locally constant functions are just constant functions hence
$S^D_{k}(x) = S^D_{l}(x)$ for $k,l
> Mm_1+L$. Similarly $S^{-D}_{k}(y) = S^{-D}_{l}(y)$ for $k,l
\geq Mm_1+L$ and all $y \in Y$. Therefore we can define $\f(x) =
S^D_m(x)$ and $\psi(y) = S^{-D}_m(y)$ where $m > Mm_1+L$. Another
obvious conclusion from Fact \ref{fi o psi = id} is that $\f \circ
\psi \equiv id_Y$ and $\psi \circ \f \equiv id_X$. If $L=0$ then
$S^D_k(x)=S^D_l(x)$ for all $k,l>0$.

\begin{fakt}\label{nierownosc f<Tf<f} There exists $m_2 \geq 0$ such that for every $\|f\| \geq m_2$ we have
$$ |f(x)| - \e(M)M\|f\|- (\e(M)+1)L \leq |Tf(\f(x))| \leq |f(x)| + \e(M)\|f\|+L.$$
If $L=0$ then we can take $m_2=0$.
\end{fakt}
\begin{proof}
Let us notice that the inequality
$$|Tf(\f(x))| \leq |f(x)| + \e(M)\|f\|+L$$
holds for $f$ such that $\|f\| \geq Mm_1+L$. Indeed it follows
from Fact \ref{nierownosc Tf<f} if we take $m=\|f\|$. In this case
$S^D_m(x)=\{\f(x)\}$. Similarly, applying Fact \ref{nierownosc
Tf<f} to $T^{-1}$ and $m=\|Tf\|$ we get
\begin{eqnarray*}
|f(x)| = |T^{-1}(Tf)(\psi(\f(x)))|& \leq & |Tf(\f(x))|+\e(M)\|Tf\|
+ L\\ & \leq & |Tf(\f(x))|+\e(M)M\|f\|+(\e(M)+1)L.
\end{eqnarray*}
It holds if $\|Tf\| \geq Mm_1+L$ which is true if $\|f\| \geq
M(Mm_1+L)+L$. In order to finish the proof put $m_2=M(Mm_1+L)+L$.
\end{proof}

\begin{fakt}
$\f$ is a continuous function from $X$ to $Y$.
\end{fakt}
\begin{proof}
Assume the contrary that there is a sequence
$(x_{\sigma})_{{\sigma} \in \Sigma}$ converging to \mbox{$x_0 \in
X$} such that $\lim_{{\sigma} \in \Sigma} \f(x_{\sigma}) \neq
\f(x_0)$. By compactness of $Y$ we can assume that $\lim_{{\sigma}
\in \Sigma} \f(x_{\sigma})=y \neq \f(x_0)$. Consider $f$ such that
$\|f\|=f(\psi(y))=m$ and $f(x_0)=0$. It can be done since $\psi(y)
\neq x_0 = \psi(\f(x_0))$.  Let us assume that $m > m_2$. Then
Fact \ref{nierownosc f<Tf<f} applied to $T$ and $x_{\sigma}$ gives
us
$$|Tf(y)| \leq |f(x_0)|+\e(M)\|f\|+L.$$
On the other hand by applying Fact \ref{nierownosc f<Tf<f} once
again to $T^{-1}$ we obtain
$$|f(\psi(y))|-\e(M)\|Tf\|-L \leq |Tf(y)|.$$
Hence combining the above inequalities and the fact that $\|Tf\|
\leq Mm+L$ we get
$$m-\e(M)m-\e(M)Mm \leq 2L + \e(M)L.$$
It is obviously impossible if $1-\e(M)-\e(M)M > 0$ provided $m$ is
large enough.
\end{proof}

\begin{proof}[Proof of Theorem \ref{tw coarse banach stone}] We are ready
to prove Theorem \ref{tw coarse banach stone}. Assume that
\mbox{$M<\sqrt{1.2}$}. Let us consider
$D=\frac{2}{M}-M-(M-1)^2<\frac{2}{M}-M$ ($D$ is positive).
Therefore spaces $X$ and $Y$ are homeomorphic if $1-\e(M)M-\e(M) >
0$ where $\e(M)=2M-1-D$. The simple analysis gives us that $D \geq
4-3M$ for $1 \leq M \leq \sqrt{1.2}$ hence $\e(M) \leq 5(M-1)$.
Using that estimation we see that the inequality $$6-5M^2 \geq
1-\e(M)M-\e(M) > 0$$ holds for $M^2 < 1.2$. In order to complete
the proof of Theorem \ref{tw coarse banach stone} let us set
\mbox{$\D=\max\{(5M-4)L, \max\{\big{|}|Tf(\varphi(x))| -
|f(x)|\big{|}\textrm{; } \|f\| \leq m_2\} \} \leq (M+1)m_2+L$}.\\
Obviously for $L=0$ also $m_2=0$ hence $\D=0$.

\end{proof}
This way we proved that whenever $d_N(C(X),C(Y)) < \frac{6}{5}$
then the compact spaces $X$ and $Y$ are homeomorphic which
improves the constant that follows from the paper of Dutrieux and
Kalton. Since $d_N(C(X),C(Y)) \leq d_{BM}(C(X),C(Y))$ (the
Banach-Mazur distance) it is well known (see \cite{Cam}) that the
constant obtained cannot be more than 2. However it is unknown if
2 is the best we can obtain in the nonlinear case (even for
bi-Lipschitz maps). We discuss this problem in the last section.
\section{Stability of coarse quasi isometries between function spaces}
In this section we estimate the distance of a coarse $(M,L)$-
quasi isometry between function spaces to the isometry between
these spaces as $M \rightarrow 1$. At this point it is worth
mentioning that for $M=1$ the problem is a classical Hyers-Ulam
problem. The result of Gevirtz \cite{Gev} shows that if there
exists a bijective coarse $(1,L)$-quasi isometry ($L$-isometry in
short) between general Banach spaces then there exists isometry of
these spaces which is no more than $5L$ from the original map (the
optimal distance from the isometry is $2L$ which was shown by
Omladi\u{c} and \u{S}emrl \cite{OmlSem}). Let us notice that these
results and the classical Banach-Stone theorem give us Theorem
\ref{tw coarse banach stone} for $M=1$. Moreover the results of
this Section easily follows. For more information see the survey
paper of Rassias \cite{Rassias}.

In all the further considerations of this Section we assume that
$1<M<\sqrt{1.2}$ and that $D=\frac{2}{M}-M-(M-1)^2$. In order to
define the isometry between function spaces which is close to the
original coarse quasi isometry we need to define a continuous
function $\lambda_m :X \mapsto \{-1,1\}$ for all $m>m_2$. The
definition requires the analysis of the proof of Fact \ref{fakt ze
nonempty}. More precisely we show that there exists such a
sequence $\{f^m_{\sigma}\}_{\sigma \in \Sigma} \in P^X_m(x)$ for
every $m>m_2 \geq m_0$ and $x \in X$ that shows that $S^D_m(x)
\neq \emptyset$. From now we assume that for a given $x$ and $m$
we have chosen such a sequence which is going to be denoted by
$\{f^m_{\sigma}\}_{\sigma \in \Sigma}$ and the corresponding
sequence $\{y^m_{\sigma}\}_{\sigma \in \Sigma}$ such that
$y^m_{\sigma} \rightarrow y \in S^D_m(x)=\{\f(x)\}$. We define
$\lambda_m(x) = {\rm sign} f^m_{\sigma}(x)$. Let us notice that
the function is well defined since from the construction of
$\{f^m_{\sigma}\}_{\sigma \in \Sigma}$ it follows that the value
of ${\rm sign} f^m_{\sigma}(x)$ does not depend on the choice of
$\sigma \in \Sigma$. It also follows from the construction that
$\forall \sigma \in \Sigma$ $Tf^m_{\sigma}(y^m_{\sigma})>0$ and
$T(-f^m_{\sigma})(y^m_{\sigma})<0$. Hence the sign of $T\lambda
f^m_{\sigma}(y^m_{\sigma})$ is the same as the sign of
$\lambda_m(x)\lambda f^m_{\sigma}(x)$ where $\lambda \in
\{-1,+1\}$. The following fact shows that the similar property
holds for all functions $f \in C(X)$ not only for $f^m_{\sigma}$.
\begin{fakt}\label{fakt o znakach}
Assume that $|f(x)| > 10(M-1)\|f\|$ and let $\|f\| = m$. Then for
sufficiently large $m>m_3$, (which depends on $M$ and $L$ only),
the sign of $Tf(\f(x))$ is the same as the sign of $\lambda_m(x)
f(x)$. If $L=0$ then $m_3=0$.
\end{fakt}
\begin{proof}
It is enough to show that if the sign of $f(x)$ is the same as the
sign of $\lambda f^m_{\sigma}(x)$ ($\lambda \in \{-1,1\}$) then
the signs of $Tf(\f(x))$ and $T(\lambda
f^m_{\sigma})(y^m_{\sigma})$ are the same as well (for some
$\sigma \in \Sigma$). Assume the contrary. Then for all $\sigma
\in \Sigma$ we have $|Tf(\f(x))-T(\lambda
f^m_{\sigma})(y^m_{\sigma})| = |Tf(\f(x))|+|T(\lambda
f^m_{\sigma})(y^m_{\sigma})|$. Therefore
\begin{eqnarray*}
|Tf(\f(x))| +Dm &\leq& \lif |Tf(y^m_{\sigma})-T(\lambda f^m_{\sigma})(y^m_{\sigma})|\\
&\leq& \lif M\|f-\lambda f^m_{\sigma}\|+L \leq Mm+L.
\end{eqnarray*}
Hence $|Tf(\f(x))| \leq (M-D)m+L$. Because $D \geq -3M+4$ for $M <
\sqrt{1.2}$ then $|Tf(\f(x))| \leq 4(M-1)m+L$. By Theorem \ref{tw
coarse banach stone} we obtain
$$|f(x)|- 5M(M-1)m -\D \leq |Tf(\f(x))| \leq 4(M-1)m+L.$$
Hence $10(M-1) \|f\| < |f(x)| \leq  (5M+4)(M-1)m + L+\D$. After
rearranging we get $(6-5M)(M-1)m <L+\D$. Obviously if $1<M <
\sqrt{1.2}$ for $m$ large enough we get a contradiction. If $L=0$
the contradiction is obtained for all $m>0$.
\end{proof}

\begin{fakt}\label{fakt o tym ze lambdy sa takie same}
There exists such $m_4 \geq 0$ that for all $m,l > m_4$ and $x \in
X$ \mbox{$\lambda_m(x)=\lambda_l(x)$}.
\end{fakt}
\begin{proof} Let us set the sequence $\{g_\sigma\}_{\sigma \in \Sigma} \in
P^X_1(x)$. Assume that $\lambda_m(x)=-\lambda_{m+\delta}(x)$ where
$m>m_3$ and $\delta >0$. Consider $h^k_{\sigma}(a) =
\lambda_k(x)kg_{\sigma}(a)$ for all $a \in X$.
 We get that
$\|h^m_{\sigma}+h^{m+\delta}_{\sigma}\| \leq \delta$. Therefore
$$|T(h^m_{\sigma})(\f(x))-T(-h^{m+\delta}_{\sigma})(\f(x))| \leq M\delta+L.$$
Using Fact \ref{fakt o znakach} we have that ${\rm
sign}T(-h^{m+\delta}_{\sigma})(\f(x))=-{\rm
sign}T(h^m_{\sigma})(\f(x))$ hence
$|T(h^m_{\sigma})(\f(x))|+|T(-h^{m+\delta}_{\sigma})(\f(x))| \leq
M\delta+L$ and finally from Theorem \ref{tw coarse banach stone}
$$2m+\delta-5(M^2-M)(2m+\delta)-2\D \leq M\delta+L.$$
Rearranging the above we get
\begin{eqnarray}\label{nierownosc dla lambda}
(-10M^2+10M+2)m < (5M^2-4M-1)\delta+2\Delta +L.
\end{eqnarray}

We can clearly see that if $-10M^2+10M+2>0$ (this is the case if
$M < \sqrt{1.2}$) then for every $\delta_0 > 0$ there exists
$m_4>0$ (depending on $M$, $L$ and $\delta_0>0$) such that if $m >
m_4$ the inequality (\ref{nierownosc dla lambda}) does not hold.

Hence $\lambda_m(x)=\lambda_{m+\delta}(x)$ for $m>m_4$ and $\delta
< \delta_0$. However locally constant functions are constant thus
$\lambda_m(x)=\lambda_l(x)$ where $m,l > m_4$. If $L=0$ then we
can set $m_4=0$ by taking $\delta_0 \rightarrow 0$.
\end{proof}
From now on we denote by $\lambda(x)=\lambda_m(x)$ for some
$m>m_4$. The following theorem is a kind of summary of all the
above facts:
\begin{tw}\label{stability theorem}
Let $X$ and $Y$ be compact spaces and $C(X)$, $C(Y)$ Banach spaces
of continuous real valued functions on  $X$ and $Y$, respectively.
Let \mbox{$T:C(X) \mapsto C(Y)$} be a bijective coarse $(M,L)$-
quasi isometry such that $T(0)=0$. Then for \mbox{every} \mbox{$M
< \sqrt{1.2}$} there is a homeomorphism $\varphi:X \mapsto Y$, a
constant $m$ and a continuous function \mbox{$\lambda : X \mapsto
\{-1,+1\}$} such that whenever $|f(x)|
> 10(M-1)\|f\|$ then:
$$\big{|}Tf(\varphi(x)) - \lambda(x)f(x)\big{|} \leq 5(M^2-M)\|f\| + \D$$
for every $x \in X$ and $f \in C(X)$. $\D$ depends on $M$ and $L$
only and $\D=0$ whenever $L=0$.
\end{tw}
\begin{proof}
Because of Theorem \ref{tw coarse banach stone}, Fact \ref{fakt o
znakach} and Fact \ref{fakt o tym ze lambdy sa takie same} it is
enough to prove that $\lambda$ is a continuous function. Let us
consider $x_{\sigma} \rightarrow x_0$ such that
$\lambda_m(x_{\sigma}) \rightarrow 1$ and $\lambda_m(x_0)=-1$
where $m > m_4$. Then take $f$ so that $f(x_0)=m=\|f\|$. Therefore
$|Tf(\f(x_0)) + m| \leq 5(M^2-M)m+\D$ and $|Tf(\f(x_{\sigma})) -
m| \leq 5(M^2-M)m+\D$. Hence $Tf(\f(x_{\sigma})) \rightarrow
Tf(\f(x_0))$ we obtain that $m \leq 5(M^2-M)m+\D$. This is
impossible if $1>5(M^2-M)$ (in particular if $M< \sqrt{1.2}$) for
sufficiently large $m$.

\end{proof}

Now we are ready to formulate the answer to the main problem of
this section that is we estimate the distance of a coarse
$(M,L)$-quasi isometry of function spaces from the isometry as $M
\rightarrow 1$.
\begin{wn} \label{wn o stabilizacji}
Let us assume that $T$ and $\lambda$ are as in the theorem above.
Then
$$|Tf(\f(x)) - \lambda(x)f(x)\big{|} \leq 26(M-1)\|f\| + \D$$ for every $x \in X$ and $f \in C(X)$. $\D$ depends on $M$ and $L$ only and $\D=0$ whenever $L=0$.
\end{wn}
\begin{proof}
Indeed if $|f(x)| \leq 10(M-1)\|f\|$ then
$$|Tf(\f(x))| \leq |f(x)| + 5(M^2-M)\|f\| + \D \leq 16(M-1)\|f\| +
\D.$$ Hence
\begin{eqnarray}\label{nierownosc stabilizacyjna}
 |Tf(\f(x)) - \lambda(x)f(x)\big{|} \leq 26(M-1)\|f\| +
\D
\end{eqnarray}
as the value of $\lambda(x)$ does not matter if \mbox{$|f(x)| \leq
10(M-1)\|f\|$}. In the case $|f(x)| > 10(M-1)\|f\|$ the inequality
(\ref{nierownosc stabilizacyjna}) follows from Theorem
\ref{stability theorem} and the fact that for $M < \sqrt{1.2}$ we
have $5(M^2-M) \leq 26(M-1)$.
\end{proof}

Let us notice that the above results work also in the Lipschitz
case that is when $L=0$ hence also $\Delta=0$. This way we improve
the estimations obtained by Jarosz in \cite{Jarosz}. More
precisely, from Theorem 2 in \cite{Jarosz} it follows that

$$|Tf(\f(x)) - \lambda(x)f(x)\big{|} \leq c{'}(M-1) \|f\|$$

\noindent where $c'(M-1) \rightarrow 0$ as $M \rightarrow 1$ and
$\f$ is a homeomorphism - however defined differently than in our
paper. The analysis of the Jarosz's paper gives us that
$400(M-1)^{\frac{1}{10}} \leq c'(M-1)$. It is clear then that
$c'(M-1)> 26(M-1)$ for $M< \sqrt{1.2}$ which shows that Corollary
\ref{wn o stabilizacji} gives us an improvement. It is also worth
remarking that although all the above results hold for bijective
coarse quasi isometries the analogous estimation can be obtained
for the general case using results of this section and Fact
\ref{fakt o bliskiej bijekcji}. The only thing that is going to
change is the term $\D$.

\section{Final remarks}
In this paper we deal with the problem of finding the optimal
constant $C$ for which from the inequality $d_N(C(X),C(Y))<C$ it
follows the existence of a homeomorphism between compact spaces
$X$ and $Y$. The result obtained is $C=\frac{6}{5}$. At the end of
Section 2 we have mentioned that from the linear theory it follows
that $C \leq 2$. Let us formulate the following open problem:
\begin{problem}\label{problem Banach-Stone ogolny}
Is it true that for all compact spaces $X$ and $Y$ the inequality
$d_N(C(X),C(Y)) < 2$ implies that $X$ and $Y$ are homeomorphic?
\end{problem}
When investigating the proof of Theorem \ref{tw coarse banach
stone} one can notice that if $M^2<2$ (which is the case if
$d_N(C(X),C(Y)) < 2$) then there exists a positive number $D$ such
that $D < \frac{2}{M}-M$ hence from Fact \ref{fakt ze nonempty} we
obtain a multifunction $S^D_m : X \mapsto 2^Y$ for suitably chosen
$m$. In the author's opinion there is a hope that the existence of
such a multifunction $S^D_m$ may help to answer the above problem
in positive. Perhaps the first approach could be to look at some
special classes of compact spaces such as all countable compacta
which have very well understood structure. More precisely the
following problems should be easier but still interesting
\begin{problem}
What is the answer to Problem \ref{problem Banach-Stone ogolny} if
we consider only countable compacta $X$ and $Y$? What if $X$ is a
convergent sequence i.e. $\{\frac{1}{n} \; ; n \in \N \} \cup
\{0\}$?
\end{problem}

In papers of Jarosz \cite{Jarosz} and Dutrieux and Kalton
\cite{DutKal} the results are obtained not only for $C(X)$ spaces
where $X$ is compact but also for $C_0(X)$ spaces where $X$ is
locally compact and $C_0(X)$ is the space of all continuous real
valued functions vanishing at $\infty$. Although the author was
able to obtain a slight improvement on the constant
$\frac{17}{16}$ in the locally compact case he decided not to
include it in the paper since it would make the proof more
technical without any reasonable gain.

\noindent {\bf Acknowledgments}\\
I am very grateful to the anonymous referee for suggestions which
improved the clarity of this article.


\end{document}